\documentclass[12pt,letterpaper]{amsart}

\usepackage{amssymb,amscd}

\newcommand{\m}{{\mathfrak m}}
\newcommand{\p}{{\mathfrak p}}

\newcommand{\E}{{\mathcal E}}

\newcommand{\T}{{\mathcal T}}

\newcommand{\Q}{{\mathbf Q}}

\newcommand{\Z}{{\mathbf Z}}

\newcommand{\pone}{{\mathbf P}^1}

\newcommand{\vsp}{\vspace{8pt}}

\newcommand{\vtsp}{\vspace{16pt}}
\newcommand{\rarr}{{\rightarrow}}

\newtheorem*{thmm}{Theorem}
\newtheorem*{lemm}{Lemma}

\theoremstyle{remark}
\newtheorem{remark}{Remark}

\newcommand\gal{{ \mbox{\rm Gal}}}

\newcommand{\ov}[1]{{\overline{{#1}}}}

\usepackage{fullpage}

\begin{document}

\title{On torsion sections of elliptic fibrations}

\author{Siman Wong}

\address{Department of Mathematics \& Statistics, University of Massachusetts.
	Amherst, MA 01003-9305 USA}

\email{siman@math.umass.edu}

\thanks{Supported in part by NSA}




\begin{abstract}
Let $E$ be an elliptic curve over the function field $\Q(t)$.  Suppose that
for every number field $L\not=\Q$ and every element $\tau\in L$ such
that the specialization $E_\tau$ is smooth, the curve $E_\tau$ has a
non-trivial torsion point over $L$.  We show that $E$ has a non-trivial
torsion point over $\Q(t)$.  This provides evidence in support of a question
of  Graber-Harris-Mazur-Starr on rational pseudo-sections of arithmetic
surjective morphisms.
\end{abstract}

\maketitle


\section{Introduction}

Let $X\subset\mathbf P^n$ be a variety of positive dimension defined over
a field $K$.  In \cite{mazur}, Graber-Harris-Mazur-Starr investigate various
scenarios under which the geometry of $X$ and the algebra of $K$ would
guarantee that $X$ has a $K$-rational point. One of the questions they raise
is the following \cite[Question 7, p.~540]{mazur}:
\begin{quote}
Let $K$ be a number field, and let $C$ be a smooth curve of genus $\ge 1$
over the function field $K(t)$.  Suppose that for every non-trivial algebraic
extension
$L/K$ and every element $t_0\in L$ such that $C_t$ has a smooth specialization
at $t_0$, the curve $C_{t_0}$ has a $L$-rational point.  Does
$C$ have a $K(t)$-rational point?
\end{quote}
In this note we provide evidence in support of this question.  Before we
state our result we first set up some notation.
 For the rest
of this note, $K$ is a number field, and $B$ is either $\pone$ or an elliptic
curve over $K$ of \textit{positive} Mordell-Weil rank.  Fix a $K$-rational,
ample
 Cartier
divisor $D$ on $B$.  Then for any finite extension $L/K$, define a height
function
$
H_{B, L}
$
on $B$ as follows.  If $B=\pone$, take $H_{B, L}$ to be the multiplicative
height on $\pone(L)$ \cite[p.~174]{joe3}; otherwise, take $H_{B, L}$ to be
$
H_{B, D, L}
$,
the multiplicative height on $B(L)$ with respect to the divisor $D$.  Then
the theorems of Schanuel \cite[Thm.~B.6.2]{joe3} and N\'eron
\cite[Thm.~B.6.3]{joe3} give the asymptotic estimate
$$
N_{B, L}(x)
:=
\#\{ P\in B(L):  H_{B, L}(P) < x \}
\sim
\Bigl\{
\renewcommand{\arraystretch}{1.2}
\begin{array}{ll}
  \alpha(B, L) x^2  & \text{if $B=\pone$,}
  \\
  \beta(B, D, L) (\log x)^{\text{rank}(B/L) / 2}
  &
  \text{otherwise,}
\end{array}
\renewcommand{\arraystretch}{1}
$$
where $\alpha(B, L)$ and $\beta(B, D, L)$ are explicit, positive constants.

\vsp

\begin{thmm}
With the notation as above, let $E$ be an elliptic curve over the function
field $K(B)$.  For any point $P\in B(L)$, write $E_P$ for the specialization
of $E$ at $P$.
Suppose that for every non-trivial finite extension
$
L/K
$,
there exists a constant
$
\lambda_{E, L} > 1/2
$
such that
\begin{equation}
\#\{ P\in B(L):
  \text{$H_{B, L}(P) < x$, $E_P$ is smooth and $E_P(L)_{\text{tor}} \not = 0$}
\}
\gg_{E, L} N_{B, L}(x)^{\lambda_{E, L}}.
            \label{thm}
\end{equation}
Then $E$ has a non-trivial torsion point over $K(B)$.
\end{thmm}

\vsp

\begin{remark}
We can simplify condition (\ref{thm}) by stipulating that every smooth
specialization of $E$ at every rational point of $B(L)$ have a non-trivial
torsion point.  On the other hand, for any integer
$
D\not=0, 1, -432
$,
the torsion subgroup of
$
y^2 = x^3 + D
$
over $\Q$ is $\Z/3, \Z/2$ or trivial, depending on whether $D$ is a square,
a cube, or neither \cite[p.~34]{husemoller}.  Thus the curve
$
y^2 = x^3 + t
$
has no non-trivial torsion point over $\Q(t)$, while
$
N_{\pone, \Q}(x)^{1/2}
$
of its specializations above rational points on $\pone(\Q)$ of height $<x$ do.
In particular, the bound
$
\lambda_{E, L} > 1/2
$
is optimal.
\end{remark}

\vsp

\begin{remark}
Our argument makes crucial use of Merel's theorem on
 torsion points of elliptic curves over number fields \cite{merel}.
If we have
the analogous result for Abelian varieties, we can readily extend the Theorem
to families of Abelian varieties fibered over $B$.
\end{remark}

\vsp

\section{The absolute case}

\begin{lemm}
With the notation as above,
let $L$ be a number field containing $K$, and let $E$ be an elliptic curve
defined over the function field $L(B)$.  Then the condition {\rm(\ref{thm})}
implies that $E$ has a non-trivial torsion point defined over $L(t)$.
\end{lemm}

\begin{proof}
Fix an elliptic fibration $\pi: \E\rarr\pone$ over $L$ with generic fiber
$E$.   Then every torsion point $T$ of $E(\ov{L(B)})$ corresponds to an
$L$-rational, $L$-irreducible torsion multisection $\T$ of $\pi$, such that
$T$ and $T'$ correspond to the same $\T$ if and only if $T$ and $T'$ fall
into the same
$
\gal(\ov{L(B)}/L(B))$-orbit.
In
particular, the degree of $\pi|_\T$ is precisely the cardinality of
this Galois orbit.

\vsp

Denote by
$
\Delta_\pi\subset B(\ov{L})
$
the discriminant locus of $\pi$.  Then for every
$
P\in B(L)-\Delta_\pi
$,
every $L$-rational torsion point of the smooth fiber $E_P$ belongs to a
unique $L$-rational, $L$-irreducible torsion multisection $\T$.  Thanks to
Merel's theorem \cite{merel}, only finitely many such $\T$ can have
$L$-rational points above
$
B(L)-\Delta_\pi
$.
The condition ({\ref{thm}}) then implies that there exists at least one
such $\T$ with
\begin{equation}
\#\{ P\in B(L)-\Delta_\pi:
  \text{$H_{B, L}(P) < x$
    and
        $\T(L)\cap E_P(L)_{\text{tor}}\not=\emptyset$}
\}
\gg_{E} N_{B, L}(x)^{\lambda_{E, L}}.
                \label{hypo}
\end{equation}
In particular, $\T(L)$ is infinite, so $\T$ is absolutely irreducible, and
hence a (possibly singular) curve over $L$ with geometric genus
$
p_g(\T)\le 1
$.
Since $\pi|_\T$ is a morphism, properties of height functions imply that
$
H_{\T, \pi^\ast(D), L}
$,
the height function on $\T$ with respect to the Cartier divisor $\pi^\ast(D)$,
satisfies
(cf.~\cite[Thm.~B.2.5(b), Thm.~B.3.2(b), Remark B.3.2.1(b)]{joe3})
\begin{equation}
\#\{
  P\in \T(L): H_{\T, \pi^\ast(D), L}(P) < x
\}
\gg_{\T, D, \pi} N_{B, L}(x)^{\deg(\pi|_\T)\lambda_{E, L}}.
    \label{rephrase}
\end{equation}
Suppose $p_g(\T)=1$, so $\pi|_\T$ is an $L$-isogeny between the two curves
$\T, B$
of geometric genus $1$.  That means the two curves have the same Mordell-Weil
rank over $L$.  But since 
$
\lambda_{E, L}>1/2
$
and $B(L)$ is infinite, (\ref{rephrase}) would contradict N\'eron's
theorem \cite[Thm.~B.6.3]{joe3} unless 
$
\pi|_\T
$
has degree one, i.e.~unless $\T$ is an actual, non-trivial, $L$-rational
torsion section of
$\pi$, in which case we are done.

\vsp

Next, suppose $p_g(\T)=0$, so $\T$ is $L$-birational to $\pone$ since
$
\T(L)
$
is not empty,
and the existence of the non-constant morphism $\pi|_\T$ implies that
$B=\pone$.
Denote by
$
\pi': \pone\rarr \T
$
the desingularization of $\T$.  It is defined over $L$; so does
$
\psi = \pi\circ\pi': \pone\rarr \pone
$.
Then (\ref{hypo}) implies that
\begin{equation}
\#\{ P\in \T(L):
  H_{\pone, \psi^\ast(D), L}(P) < x
\}
\gg_{E} N_{\pone, \, L}(x)^{\deg(\pi|_\T) \deg(\pi') \lambda_{E, L}}.
             \label{psi}
\end{equation}
The Picard group of $\pone$ is trivial, so properties of height functions
\cite[Thm.~B.3.2(d)]{joe3} imply that
$
H_{\pone, \psi^\ast(D), L}
$
and the standard multiplicative height $H_{\pone, L}$ differs by a positive
multiple bounded from above and from zero.  Schanuel's theorem then implies
that the left side of (\ref{psi}) is
$
\ll N_{\pone, \, L}(x)
$.
Since $\lambda_{E, L}>1/2$, this forces
$
\deg(\pi) = \deg(\pi') = 1
$,
whence $\T$ is an actual section, as desired.
\end{proof}

\vsp

\section{Proof of the Theorem}

For every non-trivial finite extension $L/K$, the Lemma furnishes a non-trivial
torsion point $T_L$ of $E$ over $L(B)$.  Suppose $E$ is \textit{not} a constant
elliptic curve over $\ov{K}$, i.e.~there does not exists an elliptic curve
$E_0$ defined over $K$ such that $E$ and $E_0$ are isomorphic over
$
\ov{K}(B)
$
Then the Mordell-Weil group of $E$ over $\ov{K}(B)$ is finitely generated
\cite[Thm.~III.6.1]{joe2},
and hence
$
E(L(B))_{\text{tor}}
$
is bounded independent of the finite extension $L/K$.  In particular, we
can find two finite extensions $L_1, L_2$ of coprime  degree over $K$ such
that $T_{L_1} = T_{L_2}$.  Then this common, non-trivial torsion point is
defined over 
$
L_1(B) \cap L_2(B) = K(B)
$,
as desired.

\vsp

Now, suppose $E$ is a constant elliptic curve over $\ov{K}$.  Then
$E(\ov{K}(B))$ is no
longer finitely generated, and we need to proceed differently.  We now give
an arithmetic argument that is in fact applicable to all $E$.

\vsp

Fix an elliptic fibration 
$
\pi: \E\rarr B
$
over $L$ with generic fiber $E$.   Fix a $K$-rational non-empty affine open
set $U\subset B$ over which $\pi$ is smooth,
and denote by $R_U$ the corresponding affine coordinate ring.  For any finite
extension $L/K$ and any
$
P\in U(L)
$,
denote by
$
\hat{R}_{U, P}
$
the completion of $R_U\otimes_K L$ at the maximal $\m_P$ corresponding
to the $L$-rational point $P$, and by
$
\hat{L}_P
$
its field of fractions.  Then the formal group argument in
\cite[Prop.~VII.3.1]{joe1}, which does \textit{not}
require that 
$
\hat{R}_{U, P}
$
have a finite residue field, implies that
$
E(\hat{L}_P)_{\text{tor}}
$
injects into
$
\E_P( \hat{R}_{U, P}/ \m_P) \simeq \E_P(L)
$
under the specialization map.  Take $L/K$ to be a non-trivial extension of the
form 
$
L = K(\sqrt{q})
$
with $q$ a rational prime, apply Merel's theorem and we see that 
$
\#E(L(B))_{\text{tor}}
$
is uniformly bounded for all such $L$.  That means we can find two such
$
L_1\not= L_2
$
with $T_{L_1}=T_{L_2}$.  Since $L_1\cap L_2 = K$, we are done.
\qed

\if 3\
{
\vsp
\hrule
\vsp

Let $\pi: S\rarr X$ be a genus one fibration with a smooth generic fiber.
Denote by
$
\epsilon: \E\rarr X
$
the associated Jacobian fibration.  Let
$
\beta: B'\rarr B
$
be a finite morphism such that
$
\pi': S\times_B B'\rarr B'
$
is $K$-isomorphic to
$
\epsilon': \E\times_B B'\rarr B'
$.
This $K$-isomorphism takes the zero section $Z$ of $\epsilon'$ to a zero
section $Z'$ of $\pi'$.  Then
$
\pi: \beta(\Z')\rarr B
$
has degree $\deg(\beta)$.

\vsp

Let $\p'$ be a $K$-rational closed point of the scheme $B'$.  Then $Z'$
intersects the fiber of $\pi'$ over $\p'$ at the identity element of this
fiber.  Let
$
\p = \beta(\p')
$.
Then $\beta'(\Z')$ intersects $S\times_B \p$ at the identity (???) with
multiplicity $\deg(\beta)$  (???)
That means the absolutely irreducible closed subscheme $\beta(Z')$ intersects
every fiber
$
S\times_B \p
$
at exactly one closed point, always with multiplity $\deg(\beta)$.  Since
$
\beta(Z')
$
is reduced, (???)  that means $\deg(\beta)=1$, as desired.

$$
\begin{diagram}
  \node[2]{\E\times_B B'}
  \arrow{sw} \arrow{sse,r}{\beta'} \arrow{se,t,-}{\sim}
  \\
  \node{\E}  
  \arrow{sse,l}{\epsilon}
  \node[2]{S\times_B B'}
  \arrow{sw,r}{\beta'}  \arrow{s,r}{\pi'}
  \\
  \node[2]{S}
  \arrow{s,l}{\pi}
  \node{B}
  \arrow{sw,r}{\beta}
  \\
  \node[2]{B}
\end{diagram}
$$

}
\fi

\vtsp

\bibliographystyle{amsalpha}

\begin{thebibliography}{A}

\bibitem{mazur}
T.~Graber, J.~Harris, B.~Mazur and J.~Starr,
Arithmetic questions related to rationally connected varieties, in:
\textit{The legacy of Niels Henrik Abel : the Abel bicentennial, Oslo, 2002},
531-542.
Springer-Verlag, 2004.

\bibitem{joe3}
M.~Hindry and J.~H.~Silverman,
\textit{Diophantine geometry: an introduction}.
Springer-Verlag, 2000.


\bibitem{husemoller}
D.~H\"usemoller,
\textit{Elliptic curves}. Springer-Verlag, 1987.


\bibitem{merel}
L.~Merel,
Bornes pour la torsion des courbes elliptiques sur les corps de nombres.
\textit{Invent.~Math.} {\bf 124} (1996) 437-449.

\bibitem{joe1}
J.~H.~Silverman,
\textit{The arithmetic of elliptic curves}.
 Springer-Verlag, 1986.


\bibitem{joe2}
J.~H.~Silverman,
\textit{Advanced topics in the arithmetic of elliptic curves}.
 Springer-Verlag, 1994.

\end{thebibliography}

\vfill

\end{document}